\documentclass[a4paper,reqno,12pt]{amsart}

\usepackage{latexsym,amssymb}
\usepackage{graphics}
\usepackage{amsmath}
\usepackage{epsfig}
\usepackage{enumerate}
\usepackage{amsthm}
\usepackage{hyperref}
\usepackage[active]{srcltx}
\usepackage{pb-diagram}
\usepackage[english,francais]{babel}

\newtheorem{thm}{Th\'eor\`eme}

\newtheorem{cor}{Corollaire}     
\newtheorem{lem}{Lemme}

\newfont{\bb}{msbm10 at 12pt}

\newcommand{\hs}{\mathbf{S}^n_+}

\newcommand\Vol{{\mathop{\rm Vol}}}

\let\<\langle
\let\>\rangle

\newcommand{\beqt}{\begin{equation}}  \newcommand{\eeqt}{\end{equation}}
\newcommand{\bal}{\begin{align}}      \newcommand{\eal}{\end{align}}
\newcommand{\ba}{\begin{array}}      \newcommand{\ea}{\end{array}}
\newcommand{\bc}{\begin{center}}     \newcommand{\ec}{\end{center}}
\newcommand{\be}{\begin{enumerate}}  \newcommand{\ee}{\end{enumerate}}
\newcommand{\beq}{\begin{eqnarray}}  \newcommand{\eeq}{\end{eqnarray}}
\newcommand{\beQ}{\begin{eqnarray*}} \newcommand{\eeQ}{\end{eqnarray*}}
\newcommand{\bi}{\begin{itemize}}    \newcommand{\ei}{\end{itemize}}
\newcommand{\bt}{\begin{tabular}}    \newcommand{\et}{\end{tabular}}

\title{Rigidit\'e conforme des h\'emisph\`eres $\mathbf{S}^4_+$ et $\mathbf{S}^6_+$}
\author{Simon Raulot}
\address
{Laboratoire de Math\'ematiques R. Salem\\
UMR $6085$ CNRS-Universit\'e de Rouen\\
Avenue de l'Universit\'e, BP.$12$\\
Technop\^ole du Madrillet\\
$76801$ Saint-\'Etienne-du-Rouvray, France}
\email{simon.raulot@univ-rouen.fr}
\date{\today}

\begin{document}
 

\maketitle

\begin{abstract}
Soit $(M,g)$ une vari\'et\'e compacte localement conform\'ement plate \`a bord totalement ombilique de dimension quatre ou six.  
On montre que si la caract\'erisitique d'Euler-Poincar\'e de $M$ est \'egale \`a $1$ et si son invariant de Yamabe est strictement positif, alors $(M,g)$ est 
conform\'ement isom\'etrique \`a l'h\'emisph\`ere standard. En combinant cet \'enonc\'e \`a un travail de Hang-Wang \cite{sphere}, on obtient un r\'esultat
de rigidit\'e pour ces h\'emisph\`eres \`a mettre en relation avec la conjecture de Min-Oo.
\end{abstract}

\vspace{0.5cm}

\selectlanguage{english}
\begin{center}
{\bf Conformal Rigidity of the hemispheres $\mathbf{S}^4_+$ and $\mathbf{S}^6_+$}
\end{center}
\begin{abstract}
Let $(M,g)$ be a four or six dimensional compact Riemannian manifold which is locally conformally flat and assume that its boundary is totally umbilical.
In this note, we prove that if the Euler characteristic of $M$ is equal to $1$ and if its Yamabe invariant is positive, then $(M,g)$ is
conformally isometric to the standard hemisphere. As an application and using a result of Hang-Wang \cite{sphere}, we prove a rigidity result for these
hemispheres regarding the Min-Oo conjecture.
\end{abstract}


\section{Introduction}


Inspir\'e par les th\'eor\`emes de la masse positive en Relativit\'e G\'en\'erale, Min-Oo proposa la conjecture suivante:

\vspace{0.2cm}
{\bf Conjecture de Min-Oo:} {\it Soit $g$ une m\'etrique riemannienne sur l'h\'emisph\`ere $\hs$ telle que:
\begin{enumerate}
\item la courbure scalaire de $g$ est sup\'erieure ou \'egale \`a $n(n-1)$,
\item la m\'etrique induite par $g$ sur $\partial\hs$ est isom\'etrique \`a la m\'etrique ronde sur $\mathbf{S}^{n-1}$,
\item le bord est totalement g\'eod\'esique pour la m\'etrique $g$,
\end{enumerate}
alors $g$ est isom\'etrique \`a la m\'etrique ronde sur $\hs$.}

\vspace{0.2cm}
Cette conjecture a fait l'objet d'une attention particuli\`ere ces derni\`eres ann\'ees et a \'et\'e prouv\'ee dans de nombreux cas particuliers. Cependant,
dans un travail r\'ecent \cite{BrendleMarquesNeves}, Brendle, Marques et Neves construisent, pour $n\geq 3$, des m\'etriques sur l'h\'emisph\`ere 
\`a courbure scalaire strictement sup\'erieure \`a $n(n-1)$ satisfaisant les conditions $(2)$ et $(3)$. De tels exemples mettent donc en d\'efaut 
la validit\'e de l'\'enonc\'e de Min-Oo. On peut alors naturellement se demander sous quelles hypoth\`eses cet \'enonc\'e est vrai. Dans \cite{HangWang}, Hang et Wang montrent qu'en rempla\c cant 
l'hypoth\`ese $(1)$ sur la courbure scalaire par son analogue sur la courbure de Ricci, la conjecture est v\'erifi\'ee. Le cas de la dimension $3$ est
trait\'e dans \cite{eichmair} avec des conditions suppl\'ementaires sur le profil isop\'erim\'etrique de $\partial M$. Dans un autre
travail \cite{sphere}, Hang et Wang s'int\'eressent aux d\'eformations conformes des m\'etriques sur l'h\'emisph\`ere standard et ils prouvent en particulier 
la conjecture pour ce type de transformations. Plus pr\'ecis\'ement, ils montrent:
\begin{thm}\label{hemisphereconf}{\rm (\cite{sphere})}
Soit $\overline{g}=e^{2u}g_{st}$ une m\'etrique dans la classe conforme de la m\'etrique ronde $g_{st}$ \`a courbure scalaire $R_{\overline{g}}\geq n(n-1)$. Si $\overline{g}=g_{st}$ le long 
de $\partial\mathbf{S}^n_+$, alors $\overline{g}=g_{st}$ sur $\mathbf{S}^n_+$.
\end{thm}

\noindent Ce r\'esultat est une des motivations de cette note. La seconde motivation est donn\'ee par la validit\'e de la conjecture de Min-Oo en dimension $2$.
Elle d\'ecoule, dans ce cas, d'un r\'esultat de Toponogov \cite{toponogov} qu'on peut \'enoncer de la mani\`ere suivante (voir aussi \cite{HangWang}):
\begin{thm}{\rm (\cite{toponogov})}
 Soit $(M^2,g)$ une surface ferm\'ee dont la courbure de Gauss satisfait $K\geq 1$. Alors toute g\'eod\'esique ferm\'ee simple de $M$ est de longueur au plus 
$2\pi$. De plus, s'il en existe une dont la longueur vaut $2\pi$ alors $M$ est isom\'etrique \`a la sph\`ere standard.
\end{thm}

Puisque toute surface est localement conform\'ement plate, il semble naturelle de se demander si la conjecture de Min-Oo est v\'erifi\'ee 
pour les vari\'et\'es localement conform\'ement plates de dimension $n\geq 3$. On montre qu'un tel 
r\'esultat est v\'erifi\'e pour $n=4$ et $n=6$. Pour cela, on s'inspire des travaux de Hebey-Vaugon \cite{hebvaug} et Gursky \cite{gurs} afin d'obtenir tout 
d'abord le r\'esultat suivant de rigidit\'e conforme pour l'h\'emisph\`ere: si $(M,g)$ est une vari\'et\'e \`a bord localement conform\'ement plate \`a invariant de Yamabe $Y_{[g]}(M)$ strictement positif 
et \`a caract\'eristique d'Euler-Poincar\'e $\chi(M)$ strictement positive alors elle est conform\'ement isom\'etrique \`a l'h\'emisph\`ere standard
(voir les Th\'eor\`emes \ref{confball4} et \ref{confball6}). On notera que la d\'emonstration que l'on donne ici utilise la r\'esolution du probl\`eme de 
Yamabe \`a bord (voir \cite{escobar:92}) et utilise donc implicitement le th\'eor\`eme de la masse positive. En combinant nos r\'esultats au 
Th\'eor\`eme \ref{hemisphereconf}, on obtient un \'enonc\'e analogue \`a celui de Min-Oo pour les vari\'et\'es localement conform\'ement 
plates de dimension $4$ et $6$:
\begin{cor}\label{minoo46}
Soit $(M,g)$ une vari\'et\'e riemannienne compacte \`a bord de dimension $4$ ou $6$ et telle que $\chi(M)=1$. 
Supposons que le bord $(\partial M,g)$ est totalement ombilique \`a courbure moyenne positive ou nulle et isom\'etrique \`a la sph\`ere 
ronde $(\mathbf{S}^{n-1},g_{st})$. Si $(M,g)$ est localement conform\'ement plate \`a courbure 
scalaire $R\geq n(n-1)$ alors $(M,g)$ est isom\'etrique \`a l'h\'emisph\`ere standard. 
\end{cor}

\noindent Notons que notre hypoth\`ese sur la g\'eom\'etrie extrins\`eque de $\partial M$ (i.e. $\partial M$ est totalement ombilique \`a courbure moyenne 
positive ou nulle) est moins restrictive que l'hypoth\`ese $(3)$. D'autre part, combin\'ee \`a l'hypoth\`ese sur la courbure scalaire,
elle permet d'assurer que la vari\'et\'e $(M,g)$ a un invariant de Yamabe 
strictement positif permettant ainsi d'appliquer les Th\'eor\`emes \ref{confball4} et \ref{confball6}.


\section{Le cas de la dimension $4$}


On consid\`ere ici $(M^4,g)$ une vari\'et\'e riemannienne compacte de dimension $4$ \`a bord lisse $\partial M$. La formule de Chern-Gauss-Bonnet pour les 
vari\'et\'es \`a bord de dimension $4$ est donn\'ee par (voir \cite{ssc} par exemple):
\begin{eqnarray}\label{chern}
32\pi^2\chi(M)=\int_M|W|^2dv+\int_M\Big(\frac{R^2}{6}-2|E|^2\Big)dv+8\int_{\partial M}\mathcal{B}ds
\end{eqnarray}

\noindent o\`u $W$, $R$ et $E:=Ric-(R/4)g$ sont respectivement le tenseur de Weyl, la courbure scalaire et la partie sans trace du tenseur de Ricci de $(M,g)$
et
\begin{eqnarray*}
\mathcal{B}=\frac{1}{2}RH-Ric(N,N)H-Riem_{\gamma\alpha\gamma\beta}S^{\alpha\beta}+\frac{1}{3}H^3-H|S|^2+\frac{2}{3}\,{\rm Tr}(S^3).
\end{eqnarray*}

\noindent Ici $Riem$, $Ric$ d\'esignent le tenseur de courbure de Riemann et le tenseur de courbure de Ricci de $(M^4,g)$ et $N$, $S$, $H$ le champ normal unitaire
sortant \`a $\partial M$, la seconde forme fondamentale et la courbure moyenne de $\partial M$ dans 
$M$ (pour la m\'etrique $g$). Le terme ${\rm Tr}(S^3)$ est d\'efinie par $S_{\alpha\beta}S_{\beta\gamma}S_{\gamma\alpha}$. Dans la suite, 
on notera $Ric_u$, $R_u,...$ pour d\'esigner ces quantit\'es dans une m\'etrique $g_u$ conforme \`a $g$ dont le facteur conforme d\'epend d'une fonction $u$. 
Puisque le tenseur de Weyl est invariant par changement conforme de la m\'etrique et que $\chi(M)$ est un invariant topologique, on obtient alors l'invariance conforme 
de la quantit\'e
\begin{eqnarray}\label{Ein}
\mathcal{F}_{2}([g]):=\int_M\Big(\frac{R^2}{96}-\frac{|E|^2}{8}\Big)dv_u+\frac{1}{2}\int_{\partial M}\mathcal{B}ds.
\end{eqnarray}

Dans un premier temps, on donne une estimation de la premi\`ere valeur propre du laplacien conforme sur les vari\'et\'es \`a bord de
dimension $4$ en fonction de l'invariant $\mathcal{F}_2([g])$. Cette in\'egalit\'e est un analogue de \cite{hijraul} dans le cadre \`a bord.
\begin{thm}\label{hrbord}
Soit $(M^4,g)$ une vari\'et\'e riemannienne compacte \`a bord totalement ombilique, alors:
\begin{eqnarray}\label{inegalite2}
\lambda_1(L)^2\geq\frac{96}{\Vol(M^4,g)}\mathcal{F}_2([g])
\end{eqnarray}

\noindent o\`u $\lambda_1(L)$ d\'esigne la premi\`ere valeur propre de l'op\'erateur de Yamabe:
\begin{equation}\label{yamabe}
\left\lbrace
\begin{array}{ll}
L f_1=6\Delta f_1+Rf_1=\lambda_1(L)f_1\quad & \textrm{sur}\;M\\
B f_{1\,|\partial M}=(\frac{\partial f_1}{\partial N}+Hf_1)_{|\partial M}=0\quad & \textrm{le long de}\;\partial M.
\end{array}
\right.
\end{equation} 

\noindent De plus, on a \'egalit\'e si et seulement si $(M^4,g)$ est conform\'ement isom\'etrique \`a une vari\'et\'e d'Einstein dont le bord est totalement 
g\'eod\'esique.
\end{thm}

Bien que la premi\`ere valeur propre $\lambda_1(L)$ n'est pas invariante par changement conforme, elle contient n\'eanmoins des informations
sur la structure conforme de la m\'etrique $g$. En effet, il est montr\'e dans \cite{escobar:92} que son signe est un invariant conforme donn\'e par 
le signe de l'invariant de Yamabe $Y_{[g]}(M)$ associ\'e \`a $(M,g)$ (voir ci-apr\`es pour une d\'efinition pr\'ecise de cet invariant).
En particulier, l'in\'egalit\'e (\ref{inegalite2}) implique que si $(M^4,g)$ est une vari\'et\'e \`a bord  
totalement ombilique telle que $\mathcal{F}_2([g])>0$ alors l'invariant de Yamabe est strictement positif ou strictement n\'egatif. De mani\`ere \'equivalente,
cela signifie que $(M^4,g)$ poss\`ede, dans sa classe conforme, une m\'etrique \`a courbure scalaire strictement positive ou strictement n\'egative, les
deux situations \'etant bien s\^ur exclusives. Pour plus de d\'etails sur ce sujet, on pourra consulter la premi\`ere partie de \cite{escobar:92}.\\
  
{\it Preuve du Th\'eor\`eme \ref{hrbord}:}
La d\'emonstration de ce r\'esultat repose sur la covariance conforme de $\mathcal{F}_2([g])$ et sur le choix d'un facteur 
conforme adapt\'e. Remarquons tout d'abord que si $g_u$ est une m\'etrique conforme \`a $g$ \`a courbure moyenne nulle 
(et donc \`a bord totalement g\'eod\'esique dans la m\'etrique $g_u$), l'invariant $\mathcal{F}_2([g])$ s'\'ecrit:
\begin{eqnarray*}
\mathcal{F}_2([g]) =  \int_M\Big(\frac{R_u^2}{96}-\frac{|E_u|^2_u}{8}\Big)dv_u
\end{eqnarray*}
ce qui donne:
\begin{eqnarray}\label{Eins4}
\mathcal{F}_2([g]) \leq  \frac{1}{96}\int_MR_u^2dv_u
\end{eqnarray}

\noindent avec \'egalit\'e si et seulement si la m\'etrique $g_u$ est Einstein. Soit maintenant $f_1$ une solution du probl\`eme \`a bord (\ref{yamabe}) qu'on
peut choisir strictement positive sur $M$ et consid\'erons la m\'etrique riemannienne $g_{f_1}=f_1^2g$ sur $M$. En utilisant les transformations conformes
de la courbure scalaire et de la courbure moyenne, on obtient:
\begin{eqnarray*}
 R_{f_1}=f_1^{-3}Lf_{1}=\lambda_1(L)f_1^{-2}\quad\text{et}\quad H_{f_1}=f_{1}^{-2}Bf_{1}=0.
\end{eqnarray*}
En particulier, le bord $(\partial M, g_{f_1})$ est totalement g\'eod\'esique dans $(M,g_{f_1})$. D'autre part, puisque $dv_{f_1}=f_1^4dv$ et que $\mathcal{F}_2([g])$ est un 
invariant conforme, l'in\'egalit\'e (\ref{Eins4}) dans la m\'etrique $g_{f_1}$ donne bien l'estimation (\ref{inegalite2}). Le cas d'\'egalit\'e s'obtient facilement puisqu'on a \'egalit\'e dans 
(\ref{Eins4}) et la m\'etrique $g_{f_1}$ est donc Einstein (\`a bord totalement g\'eod\'esique).
\hfill$\square$\\

Dans \cite{escobar:92}, Escobar introduit l'invariant de Yamabe:
\begin{eqnarray*}
Y_{[g]}(M):=\underset{g_u\in[g]\,/\,H_u=0}{\inf}\Big(\frac{\int_{M}R_udv_u}{\Vol(M,g_u)^{\frac{1}{2}}}\Big)
\end{eqnarray*}

\noindent afin d'\'etudier le probl\`eme suivant: \'etant donn\'ee une vari\'et\'e riemannienne compacte \`a bord, existe-t-il une m\'etrique
conforme (\`a la m\'etrique initiale) \`a courbure scalaire constante et \`a courbure moyenne nulle? Ce probl\`eme est un analogue du probl\`eme de Yamabe
classique dans le cadre des vari\'et\'es \`a bord. Il n'est pas difficile de v\'erifier que si $Y_{[g]}(M)\geq 0$ on a: 
\begin{eqnarray*}
Y_{[g]}(M)=\underset{g_u\in[g]}{\inf}\Big(\lambda_1(L_u)\Vol(M,g_u)^{\frac{1}{2}}\Big),
\end{eqnarray*}

\noindent o\`u $\lambda_1(L_u)$ est la premi\`ere valeur propre (\ref{yamabe}) du laplacien conforme $L_u$ pour la m\'etrique $g_u$. Ainsi en combinant 
l'estimation du Th\'eor\`eme \ref{hrbord} \`a cette caract\'erisation, on obtient:
\begin{eqnarray}\label{masspo}
96\mathcal{F}_2([g])\leq Y_{[g]}(M)^2\leq Y_{[g_{st}]}(\mathbf{S}^4_+)^2=192\pi^2.
\end{eqnarray}

\noindent La derni\`ere in\'egalit\'e est due \`a Escobar et l'\'egalit\'e est atteinte si et seulement si 
$(M^4,g)$ est conform\'ement isom\'etrique \`a l'h\'emisph\`ere standard de dimension $4$ (voir le Th\'eor\`eme $4.1$ dans (\cite{escobar:92})). \`A l'aide de ces estimations, on a alors:
\begin{thm}\label{confball4}
Soit $(M^4,g)$ une vari\'et\'e riemannienne compacte orient\'ee \`a bord de dimension $4$. On suppose que $M$ est localement conform\'ement plate et 
que son bord est totalement ombilique. Alors si $\chi(M)=1$ et $Y_{[g]}(M)>0$, la vari\'et\'e $(M^4,g)$ est conform\'ement isom\'etrique \`a 
l'h\'emisph\`ere $(\mathbf{S}^4_+,g_{st})$.
\end{thm}

\noindent {\it Preuve:} 
Remarquons tout d'abord que si $(M^4,g)$ est une vari\'et\'e riemannienne compacte orient\'ee, 
localement conform\'ement plate, \`a bord totalement ombilique la formule de Chern-Gauss-Bonnet (\ref{chern}) et l'in\'egalit\'e (\ref{masspo}) donnent:
\begin{eqnarray*}
32\pi^2\chi(M)=16\mathcal{F}_2([g])\leq 32\pi^2
\end{eqnarray*}

\noindent et donc $\chi(M)\leq 1$. Si on suppose maintenant que $\chi(M)=1$ on a alors \'egalit\'e dans (\ref{masspo}) et on conclut que $(M^4,g)$ est 
conform\'ement isom\'etrique \`a l'h\'emisph\`ere rond.
\hfill$\square$


\section{Le cas de la dimension $6$}


Dans cette section, on prouve un r\'esultat analogue \`a celui de la partie pr\'ec\'edente en dimension $6$.
\begin{thm}\label{confball6}
Soit $(M^6,g)$ une vari\'et\'e riemannienne compacte orient\'ee \`a bord de dimension $6$. On suppose que $M$ est localement conform\'ement plate et 
que son bord est totalement ombilique. Alors si $\chi(M)=1$ et $Y_{[g]}(M)>0$, la vari\'et\'e $(M^6,g)$ est conform\'ement isom\'etrique \`a 
l'h\'emisph\`ere $(\mathbf{S}^6_+,g_{st})$.
\end{thm}

La d\'emonstration de ce r\'esultat repose sur un argument de Gursky \cite{gurs} et sur la r\'esolution du probl\`eme de Yamabe pour les vari\'et\'es 
localement conform\'ement plates \`a bord totalement ombilique (voir le Th\'eor\`eme $4.1$ dans (\cite{escobar:92})). Encore une fois, on utilise la formule de Chern-Gauss-Bonnet 
pour les vari\'et\'es \`a bord de dimension $6$. On ne l'\'enonce pas ici en toute g\'en\'eralit\'e mais seulement sous une forme utile pour la d\'emonstration
de notre r\'esultat (voir \cite{ssc} pour un \'enonc\'e g\'en\'eral). Un point important ici est de remarquer que si on suppose que le bord est totalement 
g\'eod\'esique, alors tous les termes de bord de la formule de Chern-Gauss-Bonnet sont nulles (voir \cite{ssc}). Donc si $(M^6,g)$ est une vari\'et\'e
riemannienne compacte localement conform\'ement plate \`a bord totalement g\'eod\'esique, on a:
\begin{eqnarray}\label{chern6}
 256\pi^3\chi(M)=\int_M {\rm Tr}(E^3)dv-\frac{2}{5}\int_MR|E|^2dv+\frac{4}{225}\int_MR^3dv.
\end{eqnarray}

D'autre part, on rappelle (voir \cite{gurs}) que si $(M^n,g)$ est 
une vari\'et\'e localement conform\'ement plate de dimension $n$, on a:
\begin{eqnarray*}
\Delta E_{ij} & = & -\frac{1}{2}\Big(\frac{n-2}{n-1}\Big)\nabla_i\nabla_jR-\frac{1}{2n}\Big(\frac{n-2}{n-1}\Big)(\Delta R)g_{ij}\\
 & & +\frac{1}{n-2}|E|^2g_{ij}-\frac{n}{n-2}g^{\alpha\beta}E_{i\alpha}E_{j\beta}-\frac{1}{n-1}RE_{ij}
\end{eqnarray*}

\noindent o\`u $E:=Ric-(R/n)g$ est la partie sans trace du tenseur de Ricci. Si on suppose maintenant que $M$ est compacte \`a bord, alors en int\'egrant l'\'egalit\'e 
pr\'ec\'edente contre $g^{ik}g^{jl}E_{kl}$, on obtient \`a l'aide de la formule de Stokes et de la deuxi\`eme identit\'e de Bianchi:
\begin{eqnarray*}
 \int_{M} {\rm Tr}(E^3) & = & \frac{1}{4}\frac{(n-2)^3}{n^2(n-1)}\int_M|\nabla R|^2dv-\frac{n-2}{n}\int_M|\nabla E|^2dv +\frac{n-2}{n(n-1)}\int_MR|E|^2dv\\
& & +\frac{n-2}{n}\int_{\partial M}\big(\frac{1}{2}\frac{\partial|E|^2}{\partial N}-E(N,\nabla R)\big)ds.
\end{eqnarray*}
 
\noindent En combinant maintenant cette identit\'e pour $n=6$ \`a la formule de Chern-Gauss-Bonnet (\ref{chern6}), on obtient:
\begin{equation}\label{cgbtg}
\begin{aligned}
384\pi^3\chi(M) & = &\frac{2}{75}\int_MR^3dv+\frac{2}{15}\int_M|\nabla R|^2dv-\frac{4}{5}\int_MR|E|^2dv\\
  & & +\int_{\partial M}\Big(\frac{1}{2}\frac{\partial|E|^2}{\partial N}-E(N,\nabla R)\Big)ds ,
\end{aligned}
\end{equation}

\noindent pour toute vari\'et\'e de dimension $6$, localement conform\'ement plate \`a bord totalement g\'eod\'esique. On prouve alors:
\begin{lem}\label{chern6constantscalar}
Si $(M^6,g)$ est une vari\'et\'e localement conform\'ement plate \`a courbure scalaire constante et \`a bord totalement g\'eod\'esique, la formule de
Chern-Gauss-Bonnet est donn\'ee par:
\begin{eqnarray*}
384\pi^3\chi(M)=\frac{2}{75}R^3{\rm Vol}(M,g)-\frac{4R}{5}\int_M|E|^2dv.
\end{eqnarray*}
\end{lem}

\noindent {\it Preuve:}
Puisque la courbure scalaire est constante, on a $\nabla R=0$ et le terme de bord dans (\ref{cgbtg}) est donn\'e par:
\begin{eqnarray*}
\frac{1}{2}\int_{\partial M}\frac{\partial|Ric|^2}{\partial N}ds=\int_{\partial M}\<\nabla_NRic,Ric\>ds.
\end{eqnarray*}

\noindent Il suffit donc de v\'erifier que $\<\nabla_NRic,Ric\>=0$. Pour cela, on se place dans un syst\`eme de coordonn\'ees de Fermi $(x_1,...,x_n)$ au 
voisinage d'un point $p\in\partial M$. Puisque le bord est totalement g\'eod\'esique, on a en $p$:
\begin{eqnarray}\label{calculricci}
\<\nabla_NRic,Ric\>= \sum_{i,j=1}^{n-1}Ric_{ij,n}Ric_{ij}+2\sum_{i=1}^{n-1}Ric_{in,n}Ric_{in}+Ric_{nn,n}Ric_{nn}
\end{eqnarray}

\noindent o\`u $Ric_{ij,n}=(\nabla_{\partial_n}Ric)(\partial_i,\partial_j)=\partial Ric_{ij}/\partial x_n$, les directions $\partial_i:=\partial/\partial x_i$ 
\'etant tangentes \`a $\partial M$ pour $1\leq i\leq n-1$ et la direction $\partial_n:=\partial/\partial x_n$ normale \`a $\partial M$. Notons tout d'abord 
que le deuxi\`eme terme dans (\ref{calculricci}) est nul. En effet, l'\'equation de Codazzi donne pour $1\leq i,j,k\leq n-1$:
\begin{eqnarray}\label{codriem}
 Riem_{ijkn}=S_{ik,j}-S_{jk,i}=0
\end{eqnarray}

\noindent puisque le bord est totalement g\'eod\'esique et donc $Ric_{in}=0$ pour $1\leq i\leq n-1$. Pour le dernier terme de $(\ref{calculricci})$, on 
contracte deux fois la deuxi\`eme identit\'e de Bianchi et on obtient:
\begin{eqnarray*}
0=R_{,l}=2\sum_{i=1}^{n-1}Ric_{il,i}+2Ric_{nl,n}
\end{eqnarray*}

\noindent pour tout $1\leq l\leq n$ puisque $R$ est constante sur $M$. En particulier, pour $l=n$, on a $Ric_{nn,n}=0$ car $Ric_{in}=0$ sur $\partial M$. 
Examinons maintenant le premier terme. On remarque tout d'abord que pour $1\leq i,j\leq n-1$
\begin{eqnarray}\label{riem1}
 Ric_{ij,n} & = & \sum_{k=1}^nRiem_{ikjk,n}\nonumber\\
& = & -\sum_{k=1}^{n-1}Riem_{iknj,k}-\sum_{k=1}^{n-1}Riem_{ikkn,j}+Riem_{injn,n}=Riem_{injn,n}
\end{eqnarray}
o\`u on a utilis\'e la deuxi\`eme identit\'e de Bianchi et la relation (\ref{codriem}). D'autre part, puisque la vari\'et\'e $(M,g)$ est localement 
conform\'ement plate et de dimension $6$, son tenseur de Weyl s'annule et sa courbure de Riemann est donc donn\'ee par:
\begin{eqnarray*}
 Riem_{ijkl}=\frac{1}{4}\big(Ric_{ik}g_{jl}+Ric_{jl}g_{ik}-Ric_{il}g_{jk}-Ric_{jk}g_{il}\big)-\frac{R}{20}\big(g_{ik}g_{jl}-g_{il}g_{jk}\big)
\end{eqnarray*}
d'o\`u:
\begin{eqnarray*}
Riem_{injn}=\frac{1}{4}\big(Ric_{ij}g_{nn}+Ric_{nn}g_{ij}-Ric_{in}g_{jn}-Ric_{jn}g_{in}\big)-\frac{R}{20}\big(g_{ij}g_{nn}-g_{in}g_{jn}\big).
\end{eqnarray*}
En d\'erivant cette identit\'e dans la direction $\partial_n$, on obtient
\begin{eqnarray}\label{riem2}
 Riem_{injn,n}=\frac{1}{4}Ric_{ij,n}
\end{eqnarray}
o\`u on a successivement utilis\'e le fait que $1\leq i,j\leq n-1$ et que dans la m\'etrique $g$, le bord est totalement ombilique, la courbure scalaire 
est constante et $Ric_{nn,n}=0$. En comparant (\ref{riem1}) et (\ref{riem2}), on obtient $Ric_{ij,n}=0$ et donc le premier terme de (\ref{calculricci}) est 
bien nul.
\hfill$\square$\\

\noindent {\it Preuve du Th\'eor\`eme \ref{confball6}:} 
Supposons que $(M^6,g)$ ne soit pas conform\'ement isom\'etrique \`a $(\mathbf{S}^6_+,g_{st})$. Dans ce cas, par \cite{escobar:92} (le point $(ii)$ du
Th\'eor\`eme $4.1$, p.$57$), l'invariant de Yamabe de 
$(M^6,g)$ satisfait:
\begin{eqnarray*} 
0<Y_{[g]}(M)< Y_{[g_{st}]}(\mathbf{S}^6_+)=30\Big(\frac{\omega_6}{2}\Big)^{\frac{1}{3}}
\end{eqnarray*}

\noindent o\`u $\omega_6$ est le volume standard de $\mathbf{S}^6$ et il existe une m\'etrique $\overline{g}$ conforme \`a $g$ v\'erifiant 
$R_{\overline{g}}=Y_{[g]}(M)$, $H_{\overline{g}}=0$ et 
${\rm Vol}(M,\overline{g})=1$. Puisque le bord est totalement ombilique dans $g$ et que cette propri\'et\'e est invariante par 
changement conforme, le bord est totalement g\'eod\'esique pour la m\'etrique $\overline{g}$ (car $H_{\overline{g}}=0$). On peut alors appliquer 
le lemme \ref{chern6constantscalar} et comme $\chi(M)=1$ on a:
\begin{eqnarray*}
384\pi^3\leq\frac{2}{75}Y_{[g]}(M)^3<\frac{2}{75}Y_{[g_{st}]}(\mathbf{S}^6_+)^3=384\pi^3,
\end{eqnarray*}

\noindent ce qui conduit \`a une contradiction et donc $(M^6,g)$ est conform\'ement isom\'etrique \`a $(\mathbf{S}^6_+,g_{st})$.
\hfill$\square$


\section{Une remarque sur les dimensions sup\'erieures}


Dans \cite{guanlinwang} (Corollaire $1$), un r\'esultat de rigidit\'e conforme pour les vari\'et\'es localement conform\'ement plates de dimension paire est obtenu. Cependant, il
n\'ecessite des hypoth\`eses suppl\'ementaires sur la structure conforme des vari\'et\'es mises en jeu. On peut obtenir un tel \'enonc\'e dans le cadre des 
vari\'et\'es localement conform\'ement plates de dimension paire \`a bord totalement ombilique en utilisant des r\'esultats de Chen. L\`a aussi des hypoth\`eses 
suppl\'ementaires sont n\'ecessaires (voir Th\'eor\`eme $4$ dans \cite{ssc}). L'\'enonc\'e de Min-Oo peut alors \^etre v\'erifi\'e dans ce cadre \`a l'aide du 
Th\'eor\`eme \ref{hemisphereconf}. Une question int\'eressante est de savoir si ces hypoth\`eses sont r\'eellement n\'ecessaires.


\bibliographystyle{amsalpha} 
\bibliography{HemisphereBib}


\end{document}